\documentclass[12pt,twoside]{article}
\usepackage{amsmath}
\usepackage{amssymb}
\usepackage{enumerate,epsfig}
\newtheorem{lemma}{Lemma}[section]
\newtheorem{theorem}{Theorem}[section]

\newtheorem{prop}{Proposition}[section]

\newcommand{\R}{ {\mathbb R} }

\newcommand{\RR}{ {\mathbb R} }

\newcommand{\pe}{{\varphi_\varepsilon}}
\makeatletter
\@addtoreset{equation}{section}
\makeatother

\newcommand{\cqfd}{{\unskip\kern 6pt\penalty 500
\raise -2pt\hbox{\vrule\vbox to 6pt{\hrule width 6pt
\vfill\hrule}\vrule}\par}}
\newcommand{\eps}{{\varepsilon}}
\newcommand{\ind}{{\mathbb I}}

\begin{document}
\title{The evolutionary limit for models of
  populations interacting competitively with many resources} 

\date{}
\author{Nicolas Champagnat$^{1}$, Pierre-Emmanuel Jabin$^{1,2}$}

\footnotetext[1]{TOSCA project-team, INRIA Sophia Antipolis --
  M\'editerran\'ee, 2004 rte des Lucioles, BP.\ 93, 06902 Sophia
  Antipolis Cedex, France, \\
E-mail: \texttt{Nicolas.Champagnat@sophia.inria.fr}}

\footnotetext[2]{Laboratoire J.-A. Dieudonn\'e, Universit\'e de Nice --
  Sophia Antipolis, Parc Valrose, 06108 Nice Cedex 02, France, E-mail:
  \texttt{jabin@unice.fr}}
\maketitle{}

\begin{abstract}
  We consider a integro-differential nonlinear model that describes
  the evolution of a population structured 
  by a quantitative trait. The interactions between traits occur from
  competition for resources whose 
  concentrations depend on the current state of the
  population. Following the formalism of~\cite{DJMP}, we 
  study a concentration phenomenon arising in the limit of strong
  selection and small mutations. We prove that 
  the population density converges to a sum of Dirac masses
  characterized by the solution $\varphi$ of a 
  Hamilton-Jacobi equation which depends on resource concentrations
  that we fully characterize in terms of the 
  function $\varphi$.
\end{abstract}

\bigskip


\section{Introduction}
We are interested in the dynamics of a population subject to mutation
and selection driven by competition for 
resources. Each individual in the population is characterized by a
quantitative phenotypic trait $x\in\RR$ 
(for example the size of individuals, their age at maturity, or their
rate of intake of nutrients). 

We study the following equation
\begin{equation}
\partial_t u_\eps(t,x)=\frac{1}{\eps}\left(\sum_{i=1}^k
I_i^\eps(t)\,\eta_i(x)-1\right)\,u_\eps(t,x)+M_\eps(u_\eps)(t,x), \label{ueps}
\end{equation}
where $M_\eps$ is the mutation kernel
\begin{equation}
M_\eps(f)(x)=\frac{1}{\eps}\int_\R K(z)\,(f(x+\eps z)-f(x))\,dz,\label{Meps}
\end{equation}
for a $K\in C^\infty_c(\R)$ such that $\int_\RR zK(z)\,dz=0$. Among
many other ecological 
situations~\cite{diekmann-al-01}, this model is relevant for the
evolution of bacteria in a 
chemostat~\cite{diekmann-04,DJMP}. With this interpretation,
$u^{\varepsilon}(t,x)$ represents the 
concentration of bacteria with trait $x$ at time $t$, the function
$\eta_i(x)$ represents the growth rate of 
the population of trait $x$ due to the consumption of a resource whose
concentration is $I^\varepsilon_i$, and 
the term $-1$ corresponds to the decrease of the bacteria
concentration due to the constant flow out of the 
chemostat. This model extends the one proposed in~\cite{DJMP} to an
arbitrary number of resources. 

This equation has to be coupled with equations for the resources $I_i$, namely
\begin{equation}
I_i(t)=\frac{1}{1+\int_\R \eta_i(x)\,u_\eps(x)\,dx}.\label{Ii}
\end{equation}
This corresponds to an assumption of fast resources dynamics with
respect to the evolutionary dynamics. The 
resources concentrations are assumed to be at a (quasi-)equilibrium at
each time $t$, which depends on the 
current concentrations $u^\varepsilon$.

The limit $\varepsilon\rightarrow 0$ corresponds to a simultaneous
scaling of fast selection and small 
mutations. It was already considered in~\cite{DJMP}. The following
argument explains what limit behaviour for 
$u^\varepsilon$ can be expected when $\varepsilon\rightarrow
0$. Defining $\varphi_\eps$ as 
\begin{equation}
u_\eps=e^{\varphi_\eps/\eps},\quad \mbox{or}\ \varphi_\eps=\eps\,\log u_\eps,
\end{equation}
one gets the equation 
\begin{equation}
\partial_t \varphi_\eps=\sum_{i=1}^k
I_i^\eps(t)\,\eta_i(x)-1+H_\eps(\varphi_\eps),\label{phieps}
\end{equation}
where
\begin{equation}
H_\eps(f)=\int_\R
K(z)\,\left(e^{(f(x+\eps\,z)-f(x))/\eps}-1\right)\,dz.\label{Heps} 
\end{equation}
At the limit $\eps\rightarrow 0$ the Hamiltonian $H_\eps$ simply becomes
\begin{equation}
H(p)=\int_\RR K(z)\,\left(e^{p\,z}-1\right)\,dz.\label{defH}
\end{equation}
So one expects Eq.~\eqref{phieps} to lead to
\begin{equation}
\partial_t \varphi=\sum_{i=1}^k
I_i(t)\,\eta_i(x)-1+H(\partial_x\varphi),\label{phi}
\end{equation}
for some $I_i$ which are unfortunately unknown since one cannot pass
to the limit directly 
in~\eqref{Ii}. Therefore one needs to find a relation between the
$\varphi$ and the $I_i$ at the limit. Under 
quite general assumptions on the parameters (see Lemma~\ref{apriori}
below), the total population mass 
$\int_\RR u_{\varepsilon}$ is uniformly bounded over time. This
suggests that $\max_{x\in\RR}\varphi(t,x)=0$ 
should hold true for all $t\geq 0$. Together with~(\ref{phi}), this
gives a candidate for the limit dynamics 
as a solution to a Hamilton-Jacobi equation with Hamiltonian $H$ and
with unknown Lagrange multipliers $I_i$, 
subject to a maximum constraint. The limit population is then composed
at time $t$ of Dirac masses located at 
the maxima of $\varphi(t,\cdot)$.

This heuristics was justified in~\cite{DJMP} in the case of a single
resource (and when the resources evolve 
on the same time scale as the population), and the case of two
resources was only partly solved. The 
mathematical study of the convergence to the Hamilton-Jacobi equation
with maximum constraint and the study of 
the Hamilton-Jacobi equation itself have only be done in very specific
cases~\cite{DJMP,BP,CCP,PG}. In fact the main problem in this proposed
model is
that the number of unknowns (the resources) may easily be larger than
the dimension of the constraint (formally equal to the number of
points where $\varphi=0$). 

Our goal in this paper is to prove the convergence of
$\varphi_\varepsilon$ to a solution of~(\ref{phi}), 
where we give a full characterization of the functions $I_i$. 
Those are no more considered as Lagrange multipliers for a set of
constraints but are given by the solution $\varphi$ itself.
The new resulting model describes the evolution of a population as Dirac
masses and is formally well posed. 

The general problem of characterizing evolutionary dynamics as sums of
Dirac masses under biologically 
relevant parameter scalings is a key tool in \emph{adaptive
  dynamics}---a branch of biology studying the 
interplay between ecology and
evolution~\cite{HS,metz-nisbet-al-92,metz-geritz-al-96, 
dieckmann-law-96,CFBA}. The  
phenomenon of \emph{evolutionary branching}, where evolution drives an
(essentially) monotype population to 
subdivide into two (or more) distinct coexisting subpopulations, is
particularly relevant in this 
framework~\cite{metz-geritz-al-96,geritz-metz-al-97,geritz-kisdi-al-98}. When
the population state can be 
approximated by Dirac masses, this simply corresponds to the
transition from a population composed of a single 
Dirac mass to a population composed of two Dirac masses.

Several mathematical approaches have been explored to study this
phenomenon. One approach consists in studying 
the stationary behaviour of an evolutionary model involving a scaling
parameter for mutations, and then 
letting this parameter converge to 0. The stationary state has been
proved to be composed of one or several 
Dirac masses for various models (for 
deterministic PDE models, 
see~\cite{calsina-cuadrado-04,CCP,desvillettes-jabin-al-08,JR,GBV},
for Fokker-Planck PDEs corresponding to 
stochastic population genetics models, see~\cite{burger-bomze-96}, for
stochastic models, see~\cite{yu-07}, 
for game-theoretic models, see~\cite{cressman-hofbauer-05}). Closely
related to these works are the notions of 
ESS (evolutionarily stable strategies) and CSS (convergence stable
strategies)~\cite{metz-geritz-al-96,diekmann-04}, which allow one in
some cases to characterize stable 
stationary
states~\cite{calsina-cuadrado-04,
desvillettes-jabin-al-08,JR,cressman-hofbauer-05}.  

The other main approach consists in studying a simultaneous scaling of
mutation and selection, in order to 
obtain a limit dynamics where transitions from a single Dirac mass to
two Dirac masses could occur. Here 
again, deterministic and stochastic approaches have been explored. The
deterministic approach consists in 
applying the scaling of~(\ref{ueps}). The first formal results have been
obtained in~\cite{DJMP}. This was followed by several 
works on other models and on the corresponding
Hamilton-Jacobi 
PDE~\cite{CCP,PG}. For models of the type we consider here,
rigorous results (especially for the well posedness of the
Hamilton-Jacobi eq. at the limit) mainly only exist in the case with
just one resource, see \cite{BP} and \cite{BMP} (one resource but
multidimensional traits).

The stochastic approach is based on
individual-based models, which are related to 
evolutionary PDE models as those in~\cite{desvillettes-jabin-al-08,JR}
through a scaling of large 
population~\cite{CFM}. Using a simultaneous scaling of large
population and rare mutations, a stochatic limit 
process was obtained in~\cite{champagnat-06} in the case of a monotype
population (i.e.\ when the limit 
process can only be composed of a single Dirac mass), and in~\cite{CM}
when the limit population can be 
composed of finitely many Dirac masses.  

Finally note that the total population of individuals 
is typically very high, for bacteria for example. This is why even
stochastic models will usually take some limit with infinite
populations. Of course, this has some drawbacks. In particular the
population of individuals around a precise trait may turn out to be
low (even though the total population is large). As in the scaling
under consideration, one has growth or decay of order
$\exp(C/\varepsilon)$, this is in fact quite common. One of the most
important open problem would be to derive models that are both able of
dealing with very large populations and still treat correctly the
small subpopulations (keeping the stochastic effects or at least
truncating the population with less than 1 individual). 

There are already some attempts in this direction, mainly proposing 
models with
truncation, see \cite{PGau} and very recently \cite{MBPP}. 
For the moment however the
truncation is of the same order as the maximal subpopulation.
 
\bigskip

In order  to state our main result, we need some regularity and decay
assumptions on the $\eta_i$, namely 
\begin{equation}
\label{boundeta} 
\eta_i>0,\quad\exists \bar\eta\in C_0(\R),\quad \forall x,\quad
\sum_{i=1}^k\Big(|\eta_i(x)|+|\eta_i'(x)|+|\eta_i''(x)|\Big)\leq
\bar\eta(x),
\end{equation}
where $C_0(\R)$ is the set of continuous function, tending to $0$ as
$x\rightarrow \pm\infty$. 
 
In order to characterize the resources $I_i(t)$ involved
in~(\ref{phi}), we introduce a sort of metastable 
equilibrium. To this aim, we need an assumption on the number of
possible roots of the reproduction rate 
\begin{equation}
\exists \bar k\leq k,\;\forall\,I_1\dots I_k\in[0,1],\ \mbox{the
  function}\ \sum_{i=1}^k 
I_i\,\eta_i(x)-1\ \mbox{has at most}\ {\bar k}\ 
\mbox{roots}.\label{rootsnumber} 
\end{equation}
We also require an invertibility condition on the matrix $\eta_i(x_j)$
\begin{equation}
\forall\,x_1\dots x_{\bar k}\mbox{\ distinct, the\ } \bar k \mbox{\
  vectors of coordinates\ }(\eta_i(x_j))_{i=1\dots \bar k}\mbox{\ are
  free}\label{invertibility}  
\end{equation}
Then we may uniquely define the metastable measure associated with a
set $\omega$ by
\begin{prop} For any closed $\omega\subset \R$, there exists a unique
  finite nonnegative measure $\mu(\omega)$ satisfying\\ 
  i) $\mbox{supp}\,\mu\subset \omega$\\
  ii) denoting $\bar I_i(\mu)=1/(1+\int \eta_i(x)\,d\mu(x))$,
\[
\sum_{i=1}^k \bar I_i(\mu)\,\eta_i(x)-1\leq 0\ in\ \omega,
\quad \sum_{i=1}^k \bar
I_i(\mu)\,\eta_i(x)-1= 0\ on\ \mbox{supp}\,\mu. 
\]\label{metastable}
\end{prop}
Now the limiting $I_i$ are directly obtained by
\begin{equation}
I_i(t)=\bar I_i(\mu(\{\varphi(t,.)=0\})).\label{constraintlimit}
\end{equation}
We prove
\begin{theorem} Assume $K\in C^{\infty}_c(\RR)$, $\int_\RR
  zK(z)\,dz=0$, \eqref{boundeta}, 
  \eqref{rootsnumber}, \eqref{invertibility}, that the initial data
  $u_\eps(t=0)>0$ or $\pe(t=0)$ are 
  $C^2$, satisfy
  \begin{gather}
    \sup_\eps \int_\R u_\eps(t=0,x)\,dx<\infty,\quad
    \sup_\varepsilon\|\partial_x\varphi_\varepsilon(t=0,\cdot)\|_{L^\infty(\R)}
<\infty, 
    \\ 
    \inf_\eps\inf_{x\in\RR}\partial_{xx}\varphi(t=0,x)>-\infty,
  \end{gather}
  and that $\varphi_\varepsilon(t=0,\cdot)$ converges to a function
  $\varphi^0$ for the norm 
  $\|\cdot\|_{W^{1,\infty}(\R)}$.

  Then up to an extraction in $\eps$, $\pe$ converges to some
  $\varphi$ uniformly on any compact subset of 
  $[0,T]\times\R$ and in $W^{1,p}([0,T]\times K)$ for any $T>0$,
  $p<\infty$ and any compact $K$. In 
  particular, $\varphi$ is continuous. The function $I_i^\eps$
  converges to $I_i$ in $L^p([0,T])$ for any 
  $T>0$, $p<\infty$, where $I_i$ is defined from $\varphi$ as
  in~(\ref{constraintlimit}), and $I_i$ is 
  approximately right-continuous for \textbf{all} $t\geq 0$. The
  function $\varphi$ is a solution to 
  \eqref{phi} almost everywhere in $t,x$ with initial condition
  $\varphi(t=0,\cdot)=\varphi^0$. Moreover if 
  one defines $\psi=\varphi-\sum_{i=1}^k\int_0^t
  I_i(s)\,ds\,\eta_i(x)$, then $\psi$ is a viscosity solution 
  to
  \begin{equation}
    \partial_t \psi(t,x)=H\left(\partial_x \psi+\sum_{i=1}^k \int_0^t
      I_i(s)\,ds\,\eta_i(x)\right).\label{viscousphi}
  \end{equation}
  \label{theolimit}
\end{theorem}

We recall that a function $f$ on $[0,+\infty)$ is approxiamtely
  right-continuous at $t\geq 0$ if $t$ is a point 
of Lebesgue right-continuity of $f$, i.e.
$$
\lim_{s\rightarrow 0}\frac{1}{s}\int_t^{t+s}|f(\theta)-f(t)|d\theta=0.
$$

Notice that, under the assumptions of Theorem~\ref{theolimit},
$\varphi_\varepsilon(t=0,x)\rightarrow-\infty$ 
when $x\rightarrow\pm\infty$ since $\int_\RR u_\varepsilon<\infty$ and
$\varphi_\varepsilon$ is uniformly 
Lipschitz. Be also careful that we assume
$\|\varphi_\varepsilon(t=0,\cdot)-\varphi^0\|_{W^{1,\infty}(\R)}\rightarrow
0$ even though 
$\varphi_\varepsilon(t=0)$ (and thus $\varphi^0$) is not bounded.

In the proofs below, $C$ denotes a numerical constant which may change
from line to line. 
%
%
\section{Proof of Prop. \ref{metastable}}
\subsection{Uniqueness}
Assume that two measures $\mu_1$ and $\mu_2$ satisfy both points
of Prop.~\ref{metastable}. We first prove that they induce the same
ressources $\bar I_i$ and then conclude that they are equal.

\medskip

{\em 1st step: Uniqueness of the $\bar I_i$.} The argument here is
essentially an adaptation of \cite{JR}. First note that
\begin{equation}
  \label{eq:ineg-1}
  \int_\R \left(\sum_{i=1}^k \bar I_i(\mu_1)\,\eta_i(x)-1\right)\,d\mu_2+\int_\R
  \left(\sum_{i=1}^k \bar I_i(\mu_2)\,\eta_i(x)-1\right)\,d\mu_1\leq 0, 
\end{equation}
since $\mu_1$ and $\mu_2$ are non negative and by the point $ii$,
$\sum_{i=1}^k \bar I_i(\mu_j)\,\eta_i(x)-1$ is non positive on
$\omega$ for $j=1,2$.  

On the other hand since $\sum_{i=1}^k \bar I_i(\mu_j)\,\eta_i-1$ vanishes on the
support of $\mu_j$, one has for instance
\[\begin{split}
\int_\R \left(\sum_{i=1}^k \bar I_i(\mu_1)\,\eta_i-1\right)\,d\mu_2&=
\int_\R \left(\sum_{i=1}^k (\bar I_i(\mu_1)-\bar I_i(\mu_2))\,\eta_i
\right)\,d\mu_2\\
&=\sum_{i=1}^k (\bar I_i(\mu_1)-\bar I_i(\mu_2))\;\int_\R
\eta_i\,d\mu_2\\
&=\sum_{i=1}^k (\bar I_i(\mu_1)-\bar I_i(\mu_2))\;(1/\bar I_i(\mu_2)-1),
\end{split}\]
 by the definition of $\bar I_i(\mu_2)$.

Since one has
\[
\sum_{i=1}^k (\bar I_i(\mu_1)-\bar I_i(\mu_2))\;(1/\bar
I_i(\mu_2)-1/\bar I_i(\mu_1))=\sum_{i=1}^k \frac{(\bar I_i(\mu_1)-\bar
  I_i(\mu_2))^2}{\bar I_i(\mu_1)\,\bar I_i(\mu_2)}\geq 0.
\]
one deduces from~(\ref{eq:ineg-1}) that
\begin{equation}
\bar I_i(\mu_1)=\bar I_i(\mu_2),\quad i=1\dots k.\label{uniqueIi}
\end{equation}

\medskip

{\em 2nd step: Uniqueness of $\mu$.} It is not possible to deduce that
$\mu_1=\mu_2$ directly from \eqref{uniqueIi}. This degeneracy (the
possibility of having several equilibrium measures, all corresponding
to the same environment) is the reason why we require additional
assumptions on the $\eta_i$.

First of all by Assumption \eqref{rootsnumber}, point $i$ and thanks
to \eqref{uniqueIi}, we know that $\mu_1$ 
and $\mu_2$ are supported on at most $\bar k$ 
points $\{x_1,\ldots,x_{\bar k}\}$,
which are the roots of $\sum_i 
\bar{I}_i(\mu_1)\eta_i(x)-1=\sum_i\bar{I}_i(\mu_2)\eta_i(x)-1$. Therefore
one may write 
\[
\mu_j=\sum_{l=1}^{\bar k} \alpha_l^j\,\delta_{x_l}.
\]
Now \eqref{uniqueIi} tells that $\int \eta_id\mu_1=\int \eta_id\mu_2$
which means that
\[
\sum_{l=1}^{\bar k} \alpha_l^1\,\eta_i(x_l)=\sum_{l=1}^{\bar k}
\alpha_l^2\,\eta_i(x_l),\quad \forall i=1\dots k.
\]
To conclude it remains to use condition \eqref{invertibility} and get
that $\alpha_l^1=\alpha_l^2$.
\subsubsection{Existence}
The basic idea to get existence is quite simple: Solve the
equation\footnote{Existence and uniqueness are 
  trivial for~(\ref{existence}), for example by Cauchy-Lipschitz
  theorem in the set of finite positive 
  measures equipped with the total variation norm.}
\begin{equation}
\partial_t \nu=\left(\sum_{i=1}^k \bar I_i(\nu)\,\eta_i(x)-1\right)\;\nu,
\label{existence}\end{equation}
and obtain the equilibrium measure $\mu$ as the limit of $\nu(t)$ as
$t\rightarrow +\infty$. 

This is done by considering the entropy
\begin{equation}
L(\nu)=\sum_{i=1}^k \log \bar I_i(\nu)+\int d\nu=-\sum_{i=1}^k
\log\Big(1+\int \eta_i\,d\nu\Big)+\int d\nu.\label{entropy}
\end{equation}
As $-\log$ is convex and $\eta_i\geq 0$, then $L$ itself is a convex
function of $\nu$. Moreover if $\nu(t)$ solves \eqref{existence}, one
has
\begin{equation}
\frac{d}{dt} L(\nu(t))=-\int \left(\sum_{i=1}^k \bar
I_i(\nu)\,\eta_i(x)-1\right)^2\;d\nu.\label{dissipation}
\end{equation}
Therefore one expects the limit of $\nu$ and the equilibrium measure
we are looking for to be the minimum of 
$L$.

Since the $\eta_i$ are bounded, one finds
\[
L(\nu)\geq -C+c\,\int d\nu,
\]
for two numerical constants $C$ and $c$. Consequently $L$ is bounded
from below on $M^1_+(\omega)$ the set of nonnegative Radon measures
on $\omega$. It also gets its infimum on a bounded part of
$M^1_+(\omega)$. As any ball of $M^1_+(\omega)$ is compact for the
weak-* topology (dual of continuous functions with compact support),
$L$ attains its infimum, or
\[
M_0=\{\nu\in M^1_+(\omega),\ L(\nu)\leq L(\nu')\ \forall \nu'\in
M^1_+(\omega)\}\neq \emptyset. 
\] 
Now take any  $\mu\in M_0$ then take $\nu$ the solution to
\eqref{existence} with $\nu(t=0)=\mu$. $L(\nu)$ is non increasing
and since it is already at a minimum initially, it is necessarily
constant. By \eqref{dissipation}, this means that
\[
\sum_{i=1}^k \bar I_i(\mu)\,\eta_i-1=0\ \mbox{on}\ \mbox{supp}\,\mu_j.
\] 
Hence $\mu$ is in fact a stationary solution to \eqref{existence} and
it satisfies point $i$ and the second part of point $ii$ of
Prop.~\ref{metastable}. Note by the way that the uniqueness argument
in fact tells that there is a unique element in $M_0$.

It only remains to check the first part of point $ii$. By
contradiction assume that there exists a point $x_0\in \omega$ s.t.
\[
\sum_{i=1}^k \bar I_i(\mu)\,\eta_i(x_0)-1>0.
\]
Let $\alpha>0$ and define $\nu_\alpha=\mu+\alpha \delta_{x_0}\in
M^1_+(\omega)$.  Now compute 
\[\begin{split}
L(\nu_\alpha)&=\int d\mu+\alpha-\sum_{i=1}^k \log\left(1+\int \eta_i
d\mu+\alpha \eta_i(x_0)\right)\\
&=\int d\mu+\alpha-\sum_{i=1}^k \left(\log\left(1+\int \eta_i
d\mu\right)+\frac{\alpha\,\eta_i(x_0)}{1+\int \eta_i
d\mu}+O(\alpha^2)\right)\\
&=L(\mu)-\alpha\,\left(\sum_{i=1}^k \bar
I_i(\mu)\,\eta_i(x_0)-1\right)+O(\alpha^2). 
\end{split}\]
Thus $L(\nu_\alpha)<L(\mu)$ for $\alpha$ small enough which is
impossible as $\mu$ is an absolute minimum of $L$. 

Consequently the first part of $ii$ is satisfied and the proof of
Prop.~\ref{metastable} complete.

\section{Proof of Theorem \ref{theolimit}}

\subsection{A priori estimates for Eq. \eqref{phieps}}
We denote by $BV_{\textup{loc}}(\R)$ the set of functions on $\R$ with
bounded variation on any compact subset 
of $\R$, by $M^1(\omega)$ the set of signed Radon measures on the
subset $\omega$ of $\R$ equipped with the 
total variation norm.

We show the following estimates on the solution to \eqref{phieps}
\begin{lemma} Let $\varphi_\eps$ be a solution to \eqref{phieps} with
  initial data $\varphi_\eps^0$ such that 
  $\int_\R e^{\varphi_\varepsilon^0(x)/\varepsilon}\,dx<\infty$,
  $\partial_x\varphi_\varepsilon^0\in 
  L^\infty(\R)$ and $\partial_{xx}\varphi_\varepsilon^0$ uniformly
  lower bounded. Then for any $T>0$ 
  \[\begin{split}
    &\|\partial_t\varphi_\eps\|_{L^\infty([0,T]\times\R)}+\|\partial_x
  \varphi_\eps\|_{L^\infty([0,T], 
      BV_{\textup{loc}}(\R)\cap L^\infty(\R))}+\|\partial_{tx}\varphi_\eps
    \|_{L^\infty([0,T],M^1)}\leq C_T,\\
    &\forall\,t\leq T,\,x\in \R,\quad \partial_{xx}\pe(t,x) \geq -C_T,\quad
    H_\eps(\pe)\geq -C_T \eps, \\
    &\forall\,t\leq T,\ \int_\R u_\eps(t,x)\,dx\leq C_T,\quad
    \pe(t,x)\leq C_T\eps\,\log 1/\eps. 
  \end{split}\] where $C_T$ only depends on the time $T$, $\int
  u_\eps(t=0)\,dx$, 
  $\|\partial_x\varphi_\varepsilon^0\|_{L^\infty(\R)}$ and
  $\inf_x\partial_{xx}\varphi_\varepsilon^0(x)$.\label{apriori}
\end{lemma}

\noindent{\bf Proof.} We start with the easy bound on the total mass.

\medskip

{\em Step 0: Bound on the total mass.} First
notice that because of \eqref{boundeta}, there exists $R>0$ s.t.
\[
\forall |x|>R,\quad \sum_{i=1}^k \eta_i(x)\leq 1/2.
\]
Let $\psi$ be a regular test function with support in $|x|>R$, taking
values in $[0,1]$ and equal to $1$ on 
$|x|>R+1$. Using the fact that $I^\varepsilon_i(t)\leq 1$, we compute
\[\begin{split}
\frac{d}{dt}\int_\R \psi(x)\,u_\eps(t,x)\,dx\leq
&-\frac{1}{2\eps}\int_\R
  \psi(x)\,u_\eps(t,x)\,dx\\
&+\frac{1}{\eps}\int_{\R^2} K(z) (\psi(x-\eps
  z)-\psi(x))\,u_\eps(t,x)\,dz\,dx\\
&\leq -\frac{1}{2\eps}\int_\R
  \psi(x)\,u_\eps(t,x)\,dx+C\int_\R u_\eps(t,x)\,dx.  
\end{split}\] 

On the other hand as each $\eta_i>0$, one has for some constant $C$
\[
I_i^\eps(t)=\frac{1}{1+\int \eta_i\,u_\eps\,dx}\leq \frac{C}{1+
  \int(1-\psi) u_\eps\,dx}.
\]
Therefore with the same kind of estimate
\[\begin{split}
\frac{d}{dt}\int_\R (1- \psi(x))\,u_\eps(t,x)\,dx\leq
&C\int_\R u_\eps(t,x)\,dx\\
+\frac{1}{\eps} \Bigg(\frac{C}{1+
  \int(1-\psi) u_\eps\,dx}&-1\Bigg)\,\int_\R (1-\psi(x))\,u_\eps(t,x)\,dx.\\
\end{split}\]
Summing the two
\[\begin{split}
\frac{d}{dt}\int_\R u_\eps(t,x)\,dx\leq&\frac{1}{\eps} \left(\frac{C}{1+
  \int(1-\psi) u_\eps\,dx}-1\right)\,\int_\R (1-
\psi(x))\,u_\eps(t,x)\,dx \\
&-\frac{1}{2\eps}\int_\R
  \psi(x)\,u_\eps(t,x)\,dx+C\int_\R u_\eps(t,x)\,dx.
\end{split}\] Since the sum of the first two terms of the r.h.s.\ is
negative if $\int u_\varepsilon$ is 
larger than a constant independent of $\varepsilon$, this shows that
$\int u_\eps(t,x)\,dx$ remains uniformly 
bounded on any finite time interval.

\medskip
 
{\em Step 1: Bound on $\partial_x \pe$.} This is a classical bound for
solutions to some Hamilton-Jacobi equations. Here we still have to
check that it remains true uniformly at the $\eps$ level. Compute
\[\begin{split}
&\partial_t \partial_x \pe=\sum_{i=1}^k I_i^\eps(t)\,\eta_i'(x)\\
&\ +\int
K(z)\,
e^{\frac{\pe(t,x+\eps z)-\pe(t,x)}{\eps}}\,\frac{\partial_x\pe(t,x+\eps
  z)-\partial_x\pe(t,x)}{\eps}\,dz. 
\end{split}\]
We first observe that, as $I_i^\eps\in [0,\ 1]$ and
$\sum_i|\eta'_i(x)|\leq \bar \eta(x)$ 
\begin{align*}
  |\partial_t\partial_x\varphi_\varepsilon| &
  \leq\bar{\eta}(x)+\frac{2}{\varepsilon}\int_\R 
  K(z)e^{|z|\,\|\partial_x\varphi_\varepsilon(t,\cdot)\|_{L^\infty(\R)}}
\|\partial_x\varphi_\varepsilon(t,\cdot)\|_{L^\infty(\R)}  
  \,dz \\ & \leq
  \frac{C}{\varepsilon}e^{C\|\partial_x\varphi_\varepsilon(t,\cdot)\|_{L^\infty(\rho)}}, 
\end{align*}
since $K$ has compact support. This entails
\[
\|\partial_x\varphi_\varepsilon(t,\cdot)\|_{L^\infty(\R)}\leq
\frac{C}{\varepsilon}
\int_0^te^{C\|\partial_x\varphi_\varepsilon(s,\cdot)\|_{L^\infty(\R)}}\,ds, 
\]
from which easily follows that $\partial_x\varphi_\varepsilon\in
L^\infty([0,t_\varepsilon],\R)$ for some 
$t_\varepsilon>0$, which may (for the moment) depend on $\varepsilon$.

Now we use the classical maximum principle. Fix $t\in[0,T]$ such that
$C_{\varepsilon,t}:=\|\partial_x\pe(t,\cdot)\|_{L^\infty(\R)}<\infty$. For
any $x\in\R$ such that 
$\partial_x\varphi_\varepsilon(t,x)>\sup_y
\partial_x\varphi_\varepsilon(t,y)-\alpha$, where the 
constant $\alpha>0$ will be specified later, we have
$$
\partial_t\partial_x\varphi_\varepsilon(t,x)\leq\bar{\eta}(x)+\int_\R K(z)
e^{|z|\,C_{t,\varepsilon}}\frac{\alpha}{\varepsilon}\,dz\leq
C\Big(1+\frac{\alpha}{\varepsilon}e^{C\,C_{t,\varepsilon}}\Big).
$$
Therefore, choosing $\alpha=\varepsilon e^{-C\,C_{t,\varepsilon}}$, we obtain
\[
\frac{d}{dt}\, \sup_x\partial_x\pe(t,x)\leq C,
\]
for a constant $C$ independent of $\varepsilon$. Using a similar
argument for the minimum, we deduce that 
$t_\varepsilon>T$ and that $\partial_x\varphi_\varepsilon$ is bounded
on $[0,T]\times\R$ by a constant 
depending only on $T$ and $\|\partial_x\varphi_\varepsilon^0\|_{L^\infty(\R)}$.

\medskip

{\em Step 2: First bound on $H_\eps(\pe)$ and bounds on 
$\partial_t\pe$ and $\pe$.} Simply note
that
\[\begin{split}
-1\leq H_\eps(\pe(t))(x)&=\int_\R
K(z)\,
e^{\frac{\pe(t,x+\eps z)-\pe(t,x)}{\eps}}\,dz-\int_\R K(z)\,dz\\
&
\leq \int
K(z)\,e^{|z|\,\|\partial_x \pe\|_{L^\infty([0,T],\R)}}\,dz\leq C.
\end{split}\]
Consequently, directly from Eq. \eqref{phieps},
\[
|\partial_t \pe|\leq \bar \eta(x)+C,
\] 
hence concluding the full Lipschitz bound on $\pe$.

To get the upper bound on $\pe$, simply note that because of the
uniform Lipschitz bound on $\pe$
\[
\pe(t,y)\geq \pe(t,x)-C_T\,|y-x|,
\]
so
\[
\int_\R u_\eps(t,y)\,dy\geq \int_\R e^{\pe(t,x)/\eps}\,e^{-C_T\,
|y-x|/\eps}
\,dy\geq C_T\,\eps\,e^{\pe(t,x)/\eps}.
\]
Hence the bound on the total mass yields that $\pe\leq C_T\,\eps\,\log
1/\eps$. 

\medskip

{\em Step 3: $BV$ bound on $\partial_x \pe$.} As for $\partial_x\pe$,
we begin with a maximum (actually, 
minimum) principle. First from \eqref{phieps}
\[\begin{split}
\partial_t \partial_{xx}\pe&\geq -\bar \eta(x)+\int_\R
K(z)\,
e^{\frac{\pe(t,x+\eps z)-\pe(t,x)}{\eps}}\,\frac{\partial_{xx}\pe(t,x+\eps
  z)-\partial_{xx}\pe(t,x)}{\eps}\,dz\\
& +\int_\R
K(z)\,
e^{\frac{\pe(t,x+\eps z)-\pe(t,x)}{\eps}}\,\frac{(\partial_{x}\pe(t,x+\eps
  z)-\partial_{x}\pe(t,x))^2}{\eps}\,dz. 
\end{split}\]
The last term is of course non negative and so with the same argument
as before, we get
\[
\frac{d}{dt} \inf_x \partial_{xx}\pe(t,x)\geq -C,
\]
where $C$ does not depend on $\varepsilon$. This proves the uniform
lower bound on $\partial_{xx} \pe$. On the 
other hand, for any measurable subset $A$ of $[x_1,x_2]$,
\begin{multline*}
  \int_{x_1}^{x_2} (\ind_{x\in A}-\ind_{x\not\in A})\partial_{xx}
  \pe(t,x)\,dx =\int_{x_1}^{x_2} \partial_{xx} 
  \pe(t,x)\,dx\\
\qquad  -2\,\int_{x_1}^{x_2} \partial_{xx} \pe(t,x)\,\ind_{x\not\in A}dx\\
  \leq \partial_x\pe(t,x_2)-\partial_x\pe(t,x_1) +C\,|x_2-x_1|\leq
  2\|\partial_x\pe\|_{L^\infty([0,T],\R)} 
  +C\,|x_2-x_1|.
\end{multline*}
This indeed shows that $\partial_{xx}\pe(t,\cdot)$ belongs to
$M^1([x_1,x_2])$ with total variation norm less 
than $2\|\partial_x\pe\|_{L^\infty([0,T],\R)}+C\,|x_2-x_1|$. Thus,
$\partial_{xx}\varphi_\varepsilon$ belongs to the space 
$L^\infty([0,\ T], M^1(\R))$, which entails $\partial_x\varphi_\varepsilon\in
L^\infty([0,T],BV_{\textup{loc}}(\R))$.

Finally, differentiating \eqref{phieps} once in $x$, one has
\[\begin{split}
  |\partial_{tx} \pe(t,x)|&\leq \bar \eta(x)+\int
  K(z)\,e^{|z|\,\|\partial_x\varphi_\varepsilon\|_{L^\infty}}
\frac{|\partial_{x}\pe(t,x+\eps 
    z)-\partial_{x}\pe(t,x)|}{\eps}\,dz\\
  &\leq \bar \eta(x)+C\int K(z)\,\int_0^z 
|\partial_{xx} \pe(t,x+\eps \theta)|\,d\theta dz.
\end{split}
\]
Integrating, by Fubini
\[
\begin{split}
\int_{x_1}^{x_2}|\partial_{tx} \pe|\,dx &\leq \int_{x_1}^{x_2}\bar
\eta(x)\,dx+C\,
\int_{x_1-\eps\rho}^{x_2+\eps\rho} |\partial_{xx}\pe|\,dx,
\end{split}\] where $\rho$ is such that the support of $K$ is included
in the ball centered at 0 of radius 
$\rho$. This ends the proof of all the bounds on the derivatives of $\pe$.

\medskip

{\em Conclusion.} It only remains to show the sharp lower bound on
$H_\eps(\pe)$. Let us write
\[
H_\eps(\pe)\geq \int_\R K(z)\,\exp\left(\int_0^1
z\,\partial_x\pe(t,x+\theta z\,\eps)\,d\theta\right)\,dz-\int_\R K(z)\,dz.
\]
The $BV$ bound on $\partial_x \pe$ shows that this function admits
right and left limits at all $x\in\R$. Let 
us denote $\partial_x\pe(t,x^+)$ the limit on the right and
$\partial_x\pe(t,x^-)$ the limit on the left. As 
$\partial_{xx}\pe$ is bounded from below, we know in addition that
\[
\forall x,\quad \partial_x\pe(t,x^+)\geq \partial_x\pe(t,x^-).
\]
By differentiating once more
\[\begin{split}
\int_0^1
z\,\partial_x\pe(t,x+\theta
z\,\eps)\,d\theta&\geq z\,\partial_x\pe(t,x^+)\\
&\qquad+\int_0^1z \int_0^1
\theta\,z\,\eps\,
\partial_{xx}\pe(t,x+\theta'\theta z\eps)\,d\theta'd\theta\\
&\geq z\,\partial_x\pe(t,x^+)-C\,\eps\,z^2,
\end{split}\]
again as $\partial_{xx} \pe$ is bounded from below. Finally
\begin{align*}
  H_\eps(\pe) & \geq \int_\R K(z) \exp(
  z\,\partial_x\pe(t,x^+)-C\,\eps\,z^2)\,dz-\int_\R K(z)\, dz \\ & \geq
  H(\partial_x\pe(t,x^+))-C\,\eps,
\end{align*}
where $H$ is defined as in \eqref{defH} and since $K$ is compactly
supported. Because we assumed that $\int_\R zK(z)\,dz=0$, we have
$H(p)\geq 0$ for any $p$, which ends the 
proof of Lemma~\ref{apriori}.
%
\subsection{Passing to the limit in $\pe$}
From the assumptions in Theorem~\ref{theolimit}, Lemma~\ref{apriori}
gives uniform bounds on $\pe$. 

Therefore up to an extraction in $\eps$ (still denoted with
$\varepsilon$), there exists a function $\varphi$ 
on $[0,T]\times\R$ such that $\partial_t\varphi\in
L^\infty([0,T]\times\R)$, $\partial_x \varphi\in 
L^\infty([0,T], BV_{\textup{loc}}\cap L^\infty(\R))$, $\partial_{tx}\varphi\in
L^\infty([0,T],M^1_{\textup{loc}}(\R))$ and $\partial_{xx}\varphi$
uniformly lower bounded on $[0,T]\times\R$, 
satisfying
\begin{equation}\begin{aligned}
\pe & \longrightarrow \varphi 
\quad\mbox{uniformly in\ }C(K)\mbox{\ for any
      compact\ }K\mbox{\ of\ }[0,T]\times\R,\\ 
    \partial_x\pe & \longrightarrow \partial_x \varphi \quad\mbox{in
      any\ }L^p_{\textup{loc}}([0,T],\R),\ p<\infty.
\end{aligned}\label{limitphieps}
\end{equation}
The first convergence follows from Arz{\'e}la-Ascoli theorem. For the
second convergence, observe that 
$\|\partial_x\varphi_\varepsilon\|_{L^\infty([0,T],BV_{\textup{loc}}(\R))}
+\|\partial_{tx}\varphi_\varepsilon\|_{L^\infty([0,T],M^1(\R))}\leq  
C_T$ implies that $\partial_x\varphi_\varepsilon$ is uniformly bounded
in $L^\infty([0,T]\times\R)\cap 
BV_{\textup{loc}}([0,T]\times\R)$. The convergence in
$L^p_{\textup{loc}}$ follows by compact embedding. We 
also have $\varphi\leq 0$ since otherwise the uniform bound on
$\int_\R u_\varepsilon(t,x)\,dx$ would be 
contradicted.

As the $I_i^\eps$ are bounded, it is possible to extract weak-*
converging subsequences (still denoted with $\eps$) to some $I_i(t)$.

Now, we write again
\[
H_\eps(\pe)=\int_\R K(z)\,\left(\exp\left(\int_0^1 z\,\partial_x
\pe(t,x+\eps\,z\,\theta)\,d\theta\right)-1\right)\,dz. 
\]
From the $L^\infty$ bound on $\partial_x\varphi_\varepsilon$ and its
strong convergence, one deduces that 
\begin{equation}
\label{eq:cv-Heps}
H_\eps(\pe)\longrightarrow H(\partial_x\varphi)\ \mbox{in}\ L^1_{loc}.
\end{equation} 
Therefore one may pass to the limit in \eqref{phieps} and obtain
\eqref{phi} (for the moment in the sense of 
distribution; the equality \emph{a.e.} will follow from the
convergence of $I_i^\varepsilon$ in $L^p([0,T])$, 
proved below).

In addition by following \cite{DJMP} or \cite{BP}, one may easily show
that $\psi(t,x)=\varphi(t,x)-\sum_{i=1}^k \int_0^t
I_i(s)\,ds\,\eta_i(x)$ is a viscosity solution  to
\eqref{viscousphi}. We refer the reader to \cite{DJMP} or \cite{BP}
for this technical part.

It remains to obtain \eqref{constraintlimit}, the approximate
right-continuity of $I_i$ for \emph{all} time 
$t$ and the convergence of $I_i^\epsilon$ to $I_i$ in $L^{p}([0,T])$
for $p<\infty$. This requires some sort 
of uniform continuity on the $I_i^\eps$ which is the object of the
rest of the proof. 
\subsection{Continuity in time for the $I_i^\eps$}
First of all note that, as suggested by the simulations
of~\cite{DJMP}, there are examples where the $I_i$ 
have jumps in time at the limit.
So we will only be able to prove their right-continuity.

This regularity in time comes from the stability of the equilibrium
defined through \eqref{constraintlimit} and
Prop.~\ref{metastable}. Therefore let us define
\[
\bar I_i(t)=\bar I_i(\mu(\{\varphi(t,.)\})),
\] 
where $\bar I_i$ and $\mu$ are given by Prop.~\ref{metastable} and
$\varphi$ is the uniform limit of $\pe$ as taken in the previous
subsection.

Our first goal is the following result.
\begin{lemma} For any fixed $s$, there exist  functions
  $\sigma_s,\;\tilde\sigma \in
  C(\R_+)$ with $\sigma_s(0)=\tilde\sigma(0)=0$ 
s.t.
\[
\int_s^t |I_i^\eps(r)-\bar I_i(s)|^2\,dr\leq
(t-s)\,\sigma_s(t-s)+\tilde \sigma(\eps).
\]
\label{conttime}
\end{lemma}
\noindent{\bf Remark.} Of course the whole point is that $\sigma_s$ and
$\tilde \sigma$ are uniform in $\eps$. It is also crucial for the
following that $\tilde\sigma$ does not depend on $s$.
\subsubsection{Proof of Lemma \ref{conttime}}
{\em Step 0: $\varphi$ has compact level sets.}\\
Observe that $\varphi_\varepsilon(t=0,x)\rightarrow-\infty$ when
$x\rightarrow\pm\infty$ since $\int_\R 
u_\varepsilon(t=0,x)\,dx<\infty$ and $\partial_x\varphi(t=0)$ is
bounded. Because of the uniform convergence 
of $\varphi_\varepsilon(t=0)$ to $\varphi^0$ on $\R$, one deduces that
$\varphi^0(x)\rightarrow-\infty$ when 
$x\rightarrow\pm\infty$.

Since $\partial_x\varphi\in L^\infty([0,T],\R)$ and $I_i(t)\in[0,1]$,
it follows from~(\ref{phi}) that 
$\partial_t\varphi\in L^\infty([0,T],\R)$ and thus
$\varphi(t,x)\rightarrow -\infty$ when 
$x\rightarrow\pm\infty$ for all $t\geq 0$.

Therefore, the set
$$
\Omega:=\{(t,x)\in[0,T]\times\R:\varphi(t,x)\geq -1\}
$$
is compact.

\medskip

{\em Step 1: One basic property of $\{\varphi=0\}$.}\\
Let us start by the following crucial observation
\begin{multline}
  \forall s,\ \exists \tau_s\in
  C(\R_+)\ \mbox{with}\ \tau_s(0)=0,\ \mbox{s.t.}\ \forall t\geq
  s,\ \\ \forall 
  x\in \{\varphi(t,.)=0\},\;\exists y\in
  \{\varphi(s,.)=0\}\ \mbox{with}\ |y-x|\leq \tau_s(t-s). 
\end{multline}
This is a sort of semi-continuity for $\{\varphi=0\}$. It is proved
very simply by contradiction. If it were 
not true, then
\[\begin{split}
&\exists s,\,\exists \tau_0>0,\,
\exists t_n\rightarrow s,\;t_n\geq s,\ \exists y_n\in
\{\varphi(t_n,.)=0\},\\
& d(y_n, \{\varphi(s,.)=0\})\geq \tau_0,
\end{split}\]
where $d(y,\omega)=\inf_{x\in\omega} |x-y|$ is the usual distance.

Since all the $y_n$ belong to the compact set $\Omega$ of Step~0, we can
extract a converging subsequence 
$y_n\rightarrow y$. As $\varphi$ is continuous, $\varphi(s,y)=0$ or
$y\in \{\varphi(s,.)=0\}$. On the other 
hand one would also have $d(y, \{\varphi(s,.)=0\})\geq \tau_0$ which
is contradictory. 

\medskip

{\em Step 2: The functional.}\\ Denote
\[
\mu_s=\mu(\{\varphi(s,.)=0\},
\]
as given by Prop.~\ref{metastable}. We look at the evolution of
\[
F_\eps(t)=\int_\R \log u_\eps(t,x)\,d\mu_s(x)=\frac{1}{\varepsilon}\int_\R \pe(t,x)\,d\mu_s(x),
\]
for $t\geq s$.  Compute
\[
\frac{d}{dt} F_\eps(t)=\frac{1}{\eps}\int_\R \left(\sum_{i=1}^k
I_i^\eps(t)\,\eta_i(x)-1\right)\,d\mu_s(x)+\frac{1}{\eps}\int_\R
H_\eps(\pe(t))\,d\mu_s.
\]
Now write
\[\begin{split}
\frac{1}{\eps}\int_\R \left(\sum_{i=1}^k
I_i^\eps(t)\,\eta_i(x)-1\right)&\,d\mu_s(x)=\frac{d}{dt}\int_\R
u_\eps(t,x)\,dx\\
-
\frac{1}{\eps}&\int_\R \left(\sum_{i=1}^k 
I_i^\eps(t)\,\eta_i(x)-1\right)\,(u_\eps(t,x)\,dx-d\mu_s(x)). 
\end{split}\]
As $\sum_{i=1}^k \bar I_i(s)\,\eta_i(x)-1$ vanishes on the support of
$\mu_s$, 
\[\begin{split}
  \frac{1}{\eps}\int_\R &\left(\sum_{i=1}^k
  I_i^\eps(t)\,\eta_i(x)-1\right)\,d\mu_s(x)=\frac{d}{dt}\int_\R 
  u_\eps(t,x)-\frac{A(t)}{\eps}\\
  &\qquad -\frac{1}{\eps}\int_\R \left(\sum_{i=1}^k (I_i^\eps(t)-\bar
    I_i(s))\,\eta_i(x)\right)\,(u_\eps(t,x)\,dx-d\mu_s(x)),
\end{split}\]
with
\[
A(t)=\int_\R \left(\sum_{i=1}^k
\bar I_i(s)\,\eta_i(x)-1\right)\,u_\eps(t,x)\,dx.
\]
Notice that
\[
\int_\R \left(\sum_{i=1}^k 
(I_i^\eps(t)-\bar
I_i(s))\,\eta_i(x)\right)\,(u_\eps(t,x)\,dx-d\mu_s(x))=-\sum_{i=1}^k 
\frac{(I_i^\eps(t)-\bar I_i(s))^2}{I_i^\eps(t)\,\bar I_i(s)}.
\]
So we deduce
\begin{equation}\begin{split}
\frac{1}{\eps}\int_s^t\sum_{i=1}^k 
\frac{(I_i^\eps(r)-\bar I_i(s))^2}{I_i^\eps(r)\,\bar I_i(s)}dr=&\int_\R
\log\frac{u_\eps(t,x)}{u_\eps(s,x)}\,d\mu_s-\int_\R 
(u_\eps(t,x)\!-\!u_\eps(s,x))\,dx \\
& +\int_s^t\frac{A(r)}{\eps}\,dr-\frac{1}{\eps}\int_s^t\int_\R
H_\eps(\pe(r))\,d\mu_s. 
\end{split}\label{conteps}\end{equation}

\medskip

{\em Step 3: Easy bounds.}\\ Lemma~\ref{apriori} tells that
\[
-H_\eps(\pe)\leq C_T\,\eps.
\]
The total mass stays bounded in time so
\[
-\int_\R 
(u_\eps(t,x)\!-\!u_\eps(s,x))\,dx\leq \int_\R 
(u_\eps(t,x)\!+\!u_\eps(s,x))\,dx\leq C.
\]
And furthermore
\[\begin{split}
\int_\R
\log\frac{u_\eps(t,x)}{u_\eps(s,x)}\,d\mu_s=&\frac{1}{\eps}\int_\R
(\pe(t,x)-\pe(s,x))\,d\mu_s\\
&\leq\frac{1}{\eps}\int_\R
(\varphi(t,x)-\varphi(s,x))\,d\mu_s
+\frac{2}{\eps}\,\|\pe-\varphi\|_{L^\infty(\Omega)},
\end{split}\] where the last bound comes from the fact that, by
Prop.~\ref{metastable}, $\mu_s$ is supported 
on $\{\varphi(s,.)=0\}\subset \Omega$, 
where $\Omega$ is defined in Step~0. Since in addition we know that $\varphi\leq
0$,
 \[\begin{split}
\int_\R
\log\frac{u_\eps(t,x)}{u_\eps(s,x)}\,d\mu_s
\leq \frac{2}{\eps}\,\|\pe-\varphi\|_{L^\infty(\Omega)}.
\end{split}\]
Consequently we deduce from \eqref{conteps} the bound
\begin{equation}
\frac{1}{\eps}\int_s^t\sum_{i=1}^k 
\frac{(I_i^\eps(r)-\bar I_i(s))^2}{I_i^\eps(r)\,\bar I_i(s)}dr\leq
C+\frac{2}{\eps}\,\|\pe-\varphi\|_{L^\infty(\Omega)}+\int_s^t \frac{A(r)}{\eps}\,dr.
\label{conteps2}\end{equation}

\medskip

{\em Step 4: Control on $A$ and the measure of $\{x,\, \pe\sim
  0\}$.}\\
For some $\alpha_\eps$ to be chosen later, decompose
\[\begin{split}
\int_s^t A(r)\,dr=&\int_s^t\int_\R \left(\sum_{i=1}^k \bar
I_i(s)\,\eta_i(x)-1\right)\,u_\eps(r,x)\ind_{\pe(r,x)\leq -\alpha_\eps} 
\,dx\,dr\\
&+\int_s^t \int_\R \left(\sum_{i=1}^k \bar
I_i(s)\,\eta_i(x)-1\right)\,u_\eps(r,x)\ind_{\pe(r,x)\geq -\alpha_\eps} 
\,dx\,dr.
\end{split}\]
For the first part, note again that by \eqref{boundeta}, there exists
$R$ s.t.
\[
\forall |x|>R,\ \sum_{i=1}^k \bar
I_i(s)\,\eta_i(x)\leq 1/2.
\] 
Therefore we may simply dominate
\[
\int_s^t \int_\R \left(\sum_{i=1}^k \bar
I_i(s)\,\eta_i(x)-1\right)\,u_\eps(r,x)\ind_{\pe(r,x)\leq -\alpha_\eps} 
\,dx\,dr\leq C\,(t-s)\,e^{-\alpha_\eps/\eps}.
\]
Concerning the second part, we constrain $1/2\geq\alpha_\eps\geq
\|\varphi-\pe\|_{L^\infty(\Omega)}$ and may therefore bound
\[\begin{split}
&\int_s^t \left(\sum_{i=1}^k \bar
I_i(s)\,\eta_i(x)-1\right)\,u_\eps(r,x)\ind_{\pe(r,x)\geq
  -\alpha_\eps}\\
&\qquad\quad\leq \int_s^t \left(\sum_{i=1}^k \bar
I_i(s)\,\eta_i(x)-1\right)\,u_\eps(r,x)\ind_{\varphi(r,x)\geq
  -2\alpha_\eps}.
\end{split}\]
Now $\sum_{i=1}^k \bar I_i(s)\,\eta_i(x)-1$ is nonpositive on
$\{\varphi(s,.)=0\}$ and so 
\[
\left(\sum_{i=1}^k \bar
I_i(s)\,\eta_i-1\right)\,\ind_{\varphi(r,.)=0}\leq C\,\sup_{x\in
  \{\varphi(r,.)=0\}}\,\inf_{y\in \{\varphi(s,.)=0\}}\,|y-x|\leq
C\,\tau_s(t-s), 
\] 
by Step~1 as the $\eta_i$ are uniformly Lipschitz. For two
sets $O_1$ and $O_2$, define in general
\[
\delta(O_1,O_2)=\sup_{x\in O_1}\,\inf_{y\in O_2} |x-y|.
\]
By the same argument, one gets
\[\begin{split}
\left(\sum_{i=1}^k \bar
I_i(s)\,\eta_i-1\right)\,\ind_{\varphi(r,.)\geq -2\alpha_\eps}\leq &C\,\tau_s(t-s)\\
&+C\,\delta(\{\varphi(r,.)\geq -2\alpha_\eps\},\{\varphi(r,.)=0\}).
\end{split}\]
Inequality \eqref{conteps} now becomes
\begin{equation}\begin{split}
\int_s^t&\sum_{i=1}^k 
\frac{(I_i^\eps(r)-\bar I_i(s))^2}{I_i^\eps(r)\,\bar I_i(s)}dr\leq
C\,\eps+2\,\|\pe-\varphi\|_{L^\infty(\Omega)}+C
\,(t-s)\,e^{-\alpha_\eps/\eps}\\
&+C\,\int_s^t \tau_s(r-s)\,dr+C\,\int_s^t \delta(\{\varphi(r,.)\geq
-2\alpha_\eps\},\{\varphi(r,.)=0\}).  
\end{split}\label{conteps3}\end{equation}

\medskip

{\em Conclusion.} Eq.~\eqref{conteps} indeed gives
Lemma~\ref{conttime} if one defines 
\[\begin{split}
&\sigma_s(t-s)=\frac{1}{t-s}\int_s^t\tau_s(r-s)\,ds,\\
&\tilde
\sigma(\eps)=C\eps+2\|\pe-\varphi\|_{L^\infty(\Omega)}+C\,T\,e^{-\alpha_\eps/\eps} 
\\
&\qquad+C\,\int_0^T \delta(\{\varphi(r,.)\geq
-2\alpha_\eps\},\{\varphi(r,.)=0\})\,dr.
\end{split}\]

Of course $\sigma_s$ is continuous and, as $\tau_s(0)=0$, then
trivially $\sigma_s(0)=0$.  Since 
$\{\varphi(r,.)\geq -2\alpha_\eps\}$ and $\{\varphi(r,.)=0\}$ are
subsets of $\Omega$, 
$\tilde{\sigma}(\varepsilon)$ is bounded for $\varepsilon\leq 1$, and
thus, in order to complete the proof of 
Lemma~\ref{conttime}, we only have to check that
$\tilde{\sigma}(\varepsilon)\rightarrow 0$ when 
$\varepsilon\rightarrow 0$ for a convenient choice of
$\alpha_\varepsilon$. If we take 
$\alpha_\varepsilon\geq\|\varphi_\varepsilon-\varphi\|_{L^\infty(\Omega)}$
converging to 0 slowly enough to have 
$\alpha_\varepsilon/\varepsilon\rightarrow+\infty$, we only have to prove that
\[
C\,\int_0^T \delta(\{\varphi(r,.)\geq
-2\alpha_\eps\},\{\varphi(r,.)=0\})\,dr\longrightarrow
0\quad\mbox{as\ }\eps\rightarrow 0.
\]
By dominated convergence it is enough that for any $r$
\[
\delta(\{\varphi(r,.)\geq
-2\alpha_\eps\},\{\varphi(r,.)=0\})\longrightarrow
0.
\]
Just as in Step~1 this is a direct consequence of the continuity of $\varphi$.
\subsection{Compactness of the $I_i^\eps$ and the obtention of
  \eqref{constraintlimit}} 
First notice that simply passing to the limit in Lemma~\ref{conttime}
\begin{lemma}$\exists
  \sigma_s,\;\in
  C(\R_+)$ with $\sigma_s(0)=0$ 
s.t. $\forall i$
\[
\int_s^t |I_i(r)-\bar I_i(s)|^2\,dr\leq
(t-s)\,\sigma_s(t-s).
\]\label{conttimelimit}
\end{lemma}
This means that at any point of Lebesgue continuity of
$I_i$,
one has $I_i=\bar I_i$. We recall that $a.e.$ point is a Lebesgue
point for $I_i$. As the $I_i$ were
defined only almost everywhere anyhow 
(they are weak-* limits), we may identify $I_i$ and $\bar I_i$. This
proves \eqref{constraintlimit} and that 
$\bar I_i$ is approximately continuous on the right for any time $t$
(and not only $a.e.\;t$). 

\medskip

Now let us prove the compactness in $L^1_{\textup{loc}}$ of each $I_i^\eps$. We
apply the usual criterion and hence wish to control
\[
\int_0^T\frac{1}{h}\int_s^{s+h} |I^\eps_i(t)-I^\eps_i(s)|\,dt\,ds.
\]
Decompose
\[\begin{split}
&\int_0^T\frac{1}{h}\int_s^{s+h} |I^\eps_i(t)-I^\eps_i(s)|\,dt\,ds\leq 
\int_0^T\frac{1}{h}\int_s^{s+h}
|I^\eps_i(t)-I_i(s)|\,dt\,ds\\
&\quad+\int_0^T\frac{1}{h}\int_s^{s+h}|I_i(t)-I^\eps_i(s)| \,dt\,ds
+\int_0^T\frac{1}{h}\int_s^{s+h} |I_i(t)-I_i(s)|\,dt\,ds.
\end{split}\]
The first and third terms are bounded directly from Lemmas
\ref{conttime} and \ref{conttimelimit}, for example by Cauchy-Lipschitz
\[\begin{split}
\int_0^T\frac{1}{h}\int_s^{s+h}
|I^\eps_i(t)-I_i(s)|\,dt\,ds&\leq \int_0^T \left(\frac{1}{h}\int_s^{s+h}
|I^\eps_i(t)-I_i(s)|^2\,dt\right)^{1/2}\,ds\\
&\leq \int_0^T
(\sigma_s(h)+\tilde\sigma(\eps)/h)^{1/2}\,ds.
\end{split}\]
The second term can be handled the same way after swapping the order
of integration
\[\begin{split}
\int_0^T\frac{1}{h}\int_s^{s+h}
|I^\eps_i(s)-I_i(t)|\,&dt\,ds
=\int_0^{T+h}\frac{1}{h}\int_{\max(0,t-h)}^t|I^\eps_i(s)-I_i(t)|\,ds\,dt \\
&\leq \int_0^{T+h} \left(\frac{1}{h}\int_{\max(0,t-h)}^{t}
|I^\eps_i(s)-I_i(t)|^2\,ds\right)^{1/2}\,dt\\
&\leq \int_0^{T+h}
(\sigma_t(h)+\tilde\sigma(\eps)/h)^{1/2}\,dt.
\end{split}\]
So finally we bound
\[\begin{split}
\int_0^T\frac{1}{h}\int_s^{s+h} |I^\eps_i(t)-I^\eps_i(s)|\,dt\,ds\leq &\
3\,\int_0^{T+h}(\sigma_s(h)+\tilde\sigma(\eps)/h)^{1/2}\,ds \\ & \leq
3\int_0^{T+h}\sqrt{\sigma_s(h)}\,ds+3(T+h)\sqrt{\tilde{\sigma}(\varepsilon)/h}.
\end{split}\] 
Since of course the functions $\sigma_s(\cdot)$ can be chosen
uniformly bounded in Lemma~\ref{conttime}, again 
by dominated convergence, this shows that $\forall \tau>0$, $\exists
h$, $\exists \eps_0(h)$ s.t. $\forall 
\eps<\eps_0(h)$
\[
\int_0^T\frac{1}{h}\int_s^{s+h} |I^\eps_i(t)-I^\eps_i(s)|\,dt\,ds\leq \tau.
\] 
This is enough to get compactness of the $I_i^\eps$ in $L^1_{\textup{loc}}$ and
then in any $L^p_{\textup{loc}}$ with $p<\infty$, which concludes the proof of
Theorem \ref{theolimit}.

\end{document}